\documentclass[12pt,a4paper]{article}

\usepackage[utf8]{inputenc}
\usepackage[english]{babel}
\usepackage{amsthm}
\usepackage{amsfonts}
\usepackage{amssymb}
\usepackage{amsmath}
\usepackage[affil-it]{authblk}
\usepackage{mathrsfs}
\usepackage[usenames]{color}
\usepackage{graphicx}
\usepackage{cite}

\newtheorem{Theorem}{Theorem}
\newtheorem{Lemma}{Lemma}

\newenvironment{Corollary}{\par\noindent{\bf Corollary. }\it}{\par\bigskip}

\DeclareMathOperator{\re}{Re}
\DeclareMathOperator{\im}{Im}

\DeclareMathOperator{\const}{const}
\DeclareMathOperator{\convex}{convex}

\newcommand{\CC}{\mathbb{C}}
\newcommand{\RR}{\mathbb{R}}
\newcommand{\NN}{\mathbb{N}}

\newcommand{\vint}[2]{\bigl|_{#1}^{#2}\bigr.}
\newcommand{\set}[2]{\bigl\{#1\,\bigl|\bigr.\,#2\bigr\}}

\makeatletter
\renewcommand{\section}{\@startsection{section}{1}%
{0pt}{3.5ex plus 1ex minus .2ex}%
{2.3ex plus.2ex}{\normalfont\large\bfseries}}
\makeatother

\title{Weil-Titchmarsh theory as applied to the singular non-sectorial Schrödinger operator. Conditions for discreteness of the spectrum and
compactness of the resolvent}
\author{Tumanov~S.\,N.}
\affil{Moscow Center of Fundamental and Applied Mathematics at Lomonosov Moscow State
University,
Moscow, Russia}

\begin{document}

\maketitle
\begin{abstract}
The spectral properties of the singular Schrödinger operator with complex-valued potential which takes
values in a wider region than the half-plane, have been little studied.
In general case, the operator is non-sectorial, and the numerical range coincides with the entire complex plane. In this situation we propose sufficient
conditions for discreteness of the spectrum and compactness of the resolvent.
\end{abstract}
\frenchspacing

\section{Introduction and main results}

Consider the Schrödinger operator in $L_2(\RR_+)$ given by the differential expression
$$
l(y)=-y''+q(x)y,\quad x\in\RR_+.
$$
Hereinafter $\RR_+=[0,+\infty)$, the potential $q=q(x)\in L_{1,loc}(\RR_+)$ is a complex-valued function, for all sufficiently
large $x\ge x_0\ge0$ takes values in the open sector
$$
\Pi_\kappa=\set{z\in\CC}{z\ne0,\ \arg z\in (-\pi+\kappa,\pi-\kappa)},\ \kappa\ge0.
$$

We introduce the {\it maximal} domain:
$$
D=\set{y\in L_2(\RR_+)}{y,y'\in AC_{loc}(\RR_+),\ l(y)\in L_2(\RR_+)},
$$
and the form of the boundary conditions: for $\alpha_0,\alpha_1\in\CC$, $|\alpha_0|+|\alpha_1|\ne0$ let
$$
U(y)=\alpha_0y(0)+\alpha_1y'(0).
$$

The operator $L_U$ will be defined in the domain $D_U=\set{y\in D}{U(y)=0}$, where it acts as follows:
$$
L_Uy=l(y).
$$

In general, such operators are non-sectorial.

There are three main approaches to studying the spectral properties of such operators: {\it perturbation theory methods, asymptotic methods,} and {\it analytic methods}.

Perturbation theory methods, e.g. \cite{AShcal, Mityagin}, consider
$L_U=L_{U1}+L_{U2}$ as the sum of the self-adjoint operator $L_{U1}y=-y''+\re q\,y$ and the perturbation $L_{U2}y=i\im q\, y$, that is subordinate
to the original operator $L_{U1}$ in some sense.
The spectral properties of the original self-adjoint operator are completely or partially preserved as a result of such a perturbation.

Asymptotic methods, e.g., \cite{Naimark2, KrejcSiegl, Iskin}, make it possible to study a wider class of operators in the sense of the growth of the $\im q$. In particular,
they allow to detect phenomena that are impossible in a self-adjoint theory, such as the complete absence of spectrum for operators of the type of the complex Airy operator
\cite{ArnalSiegl}. A significant disadvantage of this approach is the requirement for high smoothness of the potential $q$ and the limitation of its oscillatory properties.

For example, for the existence of WKB asymptotics, it is required that $q,q'\in AC_{loc}(\RR_+)$ and
$$
\frac{1}{q^{1/4}}\frac{d^2}{dx^2}\frac{1}{q^{1/4}}=\frac{5}{16}\frac{(q')^2}{q^{5/2}}-\frac{1}{4}\frac{q''}{q^{3/2}}\in L_1(\RR_+),
$$
see \cite{Olver}, which is quite cumbersome.

As it's known \cite{Olver}, solutions to the Cauchy problem for the homogeneous equation
\begin{equation}
\label{eqodnlambda}
-y''+(q(x)-\lambda)y=0,\quad x\in\RR_+
\end{equation}
depend analytically on the parameter $\lambda\in\CC$. The study of this dependence underlies analytical methods that impose minimal restrictions on the potential $q$.

The key place in this approach is occupied by the Weyl-Titchmarsh theory.

In the classical theory \cite{Titchmarsh} the function $q$ is continuous in $\RR_+$ and real, while $\lambda$ takes
non-real values. Depending on the dimension of the space
of square-integrable solutions to the homogeneous equation \eqref{eqodnlambda}, the theory highlights two alternatives: the {\it limit point} case
and the {\it limit circle} case.

In the limit point case, the resolvent $R_\lambda=(L_U-\lambda)^{-1}$ is uniquely defined for $\lambda\not\in\RR$.
In the limiting circle case, this requires additional boundary conditions.

For two linearly independent solutions of the homogeneous equation, for example, given by initial conditions:
$\varphi(0,\lambda)=\psi'(0,\lambda)=0$, $\varphi'(0,\lambda)=\psi(0,\lambda)=1$, the classical theory constructs the Weil function $\mu(\lambda)$ which is analytic
in the upper and lower half-planes, such that for all non-real $\lambda$ the function $w(x,\lambda)=\varphi(x,\lambda)+\mu(\lambda)\psi(x,\lambda)$ is
a square-integrable solution to the homogeneous equations --- the Weyl solution.

The analytic properties of $\mu(\lambda)$, considered in the entire plane, make it possible to study the spectral properties of $L_U$.

The classical Weil--Titchmarsh theory was generalized for complex-valued $q$.

For example, Sims \cite{Sims} considered potentials $q$ taking values only in the upper, or
only in the lower half-plane. The construction of the Weyl function is carried out for $\lambda$ from the opposite half-plane.

He expanded the classification of possibilities, and identified three alternatives:
\begin{enumerate}
\item {\it Limit point and one solution in $L_2$}.
There is a unique solution to \eqref{eqodnlambda} up to a constant factor, which lies in the weighted $L_2(\RR_+,\im(q-\lambda))$. It's unique
in the class of solutions from $L_2(\RR_+)$.
\label{ensims1}
\item {\it Limit point and all solutions in $L_2$}. There is only one, up to a constant factor
solution to \eqref{eqodnlambda} in $L_2(\RR_+,\im(q-\lambda))$. At the same time, all solutions to \eqref{eqodnlambda} lie in $L_2(\RR_+)$.
\item {\it Limit circle}. All solutions lie in the weighted $L_2(\RR_+,\im(q-\lambda))$ (and in $L_2(\RR_+)$).
\end{enumerate}

Lidskii \cite{Lidsky} applied this method to operators with potentials satisfying either the condition $\re q\ge\const$ or $\pm\im q\ge\const$ and
formulated sufficient conditions for the realization of the alternative \ref{ensims1}). If for $\re q\ge\const$ it is always realized,
then for potentials with unbounded from below real part, to realize \ref{ensims1}) one has to impose additional conditions on $q$,
such as of Sears \cite{Sears}.

The Weil--Titchmarsh theory was further developed in \cite{Brown}, where considered potentials with a range of values in some half-plane not containing
negative real semiaxis. The Weil function was constructed for $\lambda$ from the additional half-plane.

This work requires $\convex\set{q(x)+r}{r>0,\,x\in\RR_+}\ne\CC$, in this case
one can choose $K\in\CC$ and $\phi>0$ so that $\re (q-K)e^{i\phi} \ge0$ for all $r>0$ and $x\in\RR_+$.

In this case, the classification
of possibilities is identical to that described by Sims, except that $\re(q-K)e^{i\phi}$ appears instead of the weight $\im(q-\lambda)$.

Here we do not touch upon the generalizations of the Weil--Titchmarsh theory to the case of systems of differential equations, we note only the papers \cite{Birger,Brown2}.

In the known generalizations of the Weil--Titchmarsh theory, it is essential that the numerical range of the operator $L_U$ has a non-empty complement, and hence
to apply it in our case of a non-sectoral $L_U$, it is necessary to overcome this difficulty.

Along with the operator $L_U$, we consider $L_U^0$, defined on the domain $D_U^0$, consisting of functions $y\in D_U$ with compact support
(its own for each such function $y$), where
$L_U^0$ acts similarly: $L_U^0y=l(y)$.

It is known \cite{Naimark} that $D_U^0\subset D_U$ is dense in $L_2(\RR_+)$, the original operator $L_U$ is closed, and $L_U^0$ admits a closure.
It may turn out that its closure
$\widetilde{L_U^0}\ne L_U$ \cite{Lidsky}.

The case when $\widetilde{L_U^0}=L_U$ will be called {\it definite} is of particular interest.

The realization of the definite case
is equivalent to $L_U^*=\overline L_{\overline U}$, \cite{Naimark}, where $\overline L_{\overline U}$ is the operator defined by adjoint
boundary conditions $\overline{U}(y)=\overline{\alpha_0}y(0)+\overline{\alpha_1}y'(0)$ and the adjoint differential expression
$\overline {l}(y)=-y''+\overline{q(x)}y$, in the same way as $L_U$.

For the sectorial operators discussed in \cite{Brown}, the realization of the definite case
is equivalent to the \ref{ensims1})-st alternative in the Sims' classification \cite{Lidsky}[Lemma 6].

Hereinafter we denote by $C$ arbitrary positive constants, possibly different in related formulas.

Let us introduce the condition on the potential $q$, which is central to the present work. We say that $q$ satisfies the condition {\bf A} for $x\ge x_0$ if:
\begin{itemize}
\item the potential $q\in L_{1,loc}(\RR_+)\cap AC_{loc}([x_0,+\infty))$;
\item the image $q([x_0,+\infty))\subset\Pi_\kappa$, $\kappa\ge0$;
\item for $x\ge x_0$ and some $C>0$ the following estimate holds a.e.:
\begin{equation}
\label{eqA}
\re p(x) \ge C+ \frac{1}{2}\Bigl|\frac{p'(x)}{p(x)}\Bigr|,
\end{equation}
where $p(x)=\sqrt{q(x)}$, the branch of the root is chosen so that $\re p(x)>0$ for $x\ge x_0$.
\end{itemize}

The inequality \eqref{eqA} is a weaker condition on $q$ than the ones of Naimark \cite{Naimark2} or Ishkin \cite{Iskin}, where
in addition to the existence of WKB asymptotics it is required that
$|q|\to+\infty$, $\kappa>0$, $|q'/q^{3/2}|\to0$, as $x\to+\infty$. The latter leads
not only to \eqref{eqA}, but even to the stronger estimate: there exist $C,\varepsilon>0$ such that for all sufficiently large $x>x_1$:
\begin{equation}
\label{eqB}
\re p(x) \ge C+\Bigl(\frac{1}{2}+\varepsilon\Bigr)\Bigl|\frac{p'(x)}{p(x)}\Bigr|.
\end{equation}

Denote
$$
\rho=\rho(x)=\re p(x) - \frac{1}{2}\Bigl|\frac{p'(x)}{p(x)}\Bigr|.
$$

Let us define the set $\mathscr{N}$ of those $\lambda\in\CC$ for which $q_\lambda=q-\lambda$ satisfies
the condition {\bf A} for
$x\ge x_0(\lambda)$ (each $\lambda$ can have its own $x_0(\lambda)$).

\begin{Lemma}
\label{lmNOpen}
Let $q$ satisfies the condition {\bf A} for $x\ge x_0$. Then

The set $\mathscr{N}$ is open [and non-empty] in $\CC$.

The following conditions are equivalent: $\mathscr{N}=\CC$ and $|q(x)|\to+\infty$ as $x\to+\infty$.
\end{Lemma}

We will prove the Lemma \ref{lmNOpen} in the section ``Auxiliary results''.

Now we formulate the main results, which will be proved in the section ``Proofs of the Theorems''. The prime denotes differentiation with respect to $x$.

\begin{Theorem}
\label{thWeil}
Let $\mathscr{N}\ne\varnothing$, then the definite case is realized: $\widetilde{L_U^0}=L_U$.

Let $\Omega\subset\mathscr{N}$ be an open and connected set. There exists a function $\eta(x,\lambda)$ ($x\in\RR_+$, $\lambda\in\Omega$)
nontrivial for all $\lambda\in\Omega$
with the following properties:
\begin{itemize}
\item functions $\eta(x,\lambda)$ and $\eta'(x,\lambda)$ are continuous in both arguments $(x,\lambda)\in\RR_+\times\Omega$;
\item for each $x\in\RR_+$ the functions $\eta(x,\lambda)$ and $\eta'(x,\lambda)$ are single-valued analytic in $\Omega$ as functions of $\lambda$;
\item for any $\lambda\in\Omega$ the function $\eta(x,\lambda)\in L_2(\RR_+)$ as a function of $x$, and is the only non-trivial solution to
\eqref{eqodnlambda}
in $L_2(\RR_+)$ up to a constant factor (not depending on $x$, but possibly depending on $\lambda$);
\end{itemize}

The spectrum of $L_U$ in $\Omega$ is discrete (the eigenspaces are one-dimensional, the root subspaces are finite-dimensional), the eigenvalues are determined by the condition
$\mathscr{W}(\lambda)=0$, where
$$
\mathscr{W}(\lambda)=\alpha_0\eta(0,\lambda)+\alpha_1\eta'(0,\lambda),
$$
and the dimensions of the root subspaces are equal to the multiplicities of zeros of $\mathscr{W}(\lambda)$.

In $\Omega$ outside the spectrum, the resolvent $R_\lambda=(L_U-\lambda)^{-1}$ has the form:
$$
(R_\lambda f)(x)=\frac{\eta(x,\lambda)}{\mathscr{W}(\lambda)}\int\limits_0^x
\chi(\xi,\lambda)f(\xi)\,d\xi
+\frac{\chi(x,\lambda)}{\mathscr{W}(\lambda)}\int\limits_x^{+\infty}
\eta(\xi,\lambda)f(\xi)\,d\xi,
$$
where $\chi(x,\lambda)=\alpha_0\varphi(x,\lambda)-\alpha_1\psi(x,\lambda)$,
and $\varphi(x,\lambda)$ and $\psi(x,\lambda)$ are solutions to the Cauchy problem for the equation \eqref{eqodnlambda} with initial conditions:
$$
\varphi(0,\lambda)=0,\quad\varphi'(0,\lambda)=1,\quad
\psi(0,\lambda)=1,\quad\psi'(0,\lambda)=0.
$$
\end{Theorem}

We call the function $\eta(x,\lambda)$ the {\it Weyl solution} to the homogeneous equation \eqref{eqodnlambda}.
\bigskip

\begin{Corollary}
Let $q$ satisfies the condition {\bf A} for $x\ge x_0$, and $|q(x)|\to+\infty$ as $x\to+\infty$. Then $\sigma(L_U)$ --- the spectrum of $L_U$ --- is discrete.
\end{Corollary}

We note that under the conditions of the Corollary we assert only the discreteness of the spectrum.
A stronger condition, the compactness of the resolvent, is given by the following
Theorem:
\begin{Theorem}
\label{th02}
Let $q$ satisfies the condition {\bf A} for $x\ge x_0$, and $\rho\to+\infty$ as $x\to+\infty$. Then the resolvent $R_\lambda$ of the operator $L_U$ is a
completely-continuous operator
for $\lambda\not\in\sigma(L_U)$.
\end{Theorem}

The following Theorem can be useful in studying the completeness of the system of eigenfunctions and associated functions of the operator $L_U$
(see, for example, \cite{Tumanov}).
\begin{Theorem}
\label{th03}
Let $q\in C^1(\RR_+)$, $|q(x)|\to+\infty$ for $x\to+\infty$, let $\kappa>0$ and
for some $0<\delta<1$ for all $x\ge x_0$:
\begin{gather}
\notag
q(x)\in\Pi_\kappa,\\
\Bigl|\frac{q'(x)}{q^{3/2}(x)}\Bigr|<4\delta\tan^{3/2}\kappa\,\sin\frac{\kappa }{2}.
\label{eqth03kappa}
\end{gather}
Then
\begin{itemize}
\item the spectrum of $L_U$ is discrete,
and the resolvent $R_\lambda$ is a
completely-continuous operator
for $\lambda\not\in\sigma(L_U)$.

\item the Weyl solution to \eqref{eqodnlambda} (delivered by the Theorem \ref{thWeil}) for each $x\ge0$ is an entire function of $\lambda$. It is
a closed kernel in the sense of Levin \cite{Levin}, i.e. if for some $f\in L_2(\RR_+)$
\begin{equation}
\label{eqth03int}
\int\limits_0^{+\infty}\eta(x,\lambda)f(x)\,dx\equiv0,\mbox { for all }\lambda\in\CC,
\end{equation}
then $f\equiv0$.
\end{itemize}
\end{Theorem}

Since, with \eqref{eqA}, there is an embedding of weighted spaces $L_2(\rho)\subset L_2\subset L_2(1/|q|)$, the following Theorem is
of interest:

\begin{Theorem}
\label{th04}
Let $q$ satisfies the condition {\bf A} for $x\ge 0$,
then $R_0=L_U^{-1}$ is a bounded operator from $L_2(\RR_+,1/|q|)$
to $L_2(\RR_+,\rho)$, while $dR_0/dx$ (acting as $d[(R_0f)(x)]/dx$) is a bounded operator from $L_2(\RR_+,1 /|q|)$
to $L_2(\RR_+,\rho/|q|)$.
\end{Theorem}
\begin{Corollary}
Let $q$ satisfies the condition {\bf A} for $x\ge 0$, then for any $y\in D$ it follows that $y\in L_2(\RR_+,\rho)$,
$y'\in L_2(\RR_+,\rho/|q|)$.
\end{Corollary}
The alternative close to the \ref{ensims1})-st of Sims realizes under conditions of the Corollary --- except that
a weighted space appears with a weaker weight.

In the Theorem \ref{th04} and in the Corollary if the potential has the stronger estimate \eqref{eqB} instead of \eqref{eqA},
one can simplify the weights $\rho$, $\rho/|q|$ to $\re p$, $\re p/|q|$ due to estimates:
\begin{gather*}
\Bigl|\frac{p'}{p}\Bigr|\le\frac{1}{1/2+\varepsilon}\re p,\\
\rho\ge\re p-\frac{1/2}{1/2+\varepsilon}\re p=\frac{\varepsilon}{1/2+\varepsilon}\re p,
\end{gather*}
as a result:
$$
C\re p \le\rho\le\re p.
$$

The condition $\kappa>0$, regardless of \eqref{eqA} or \eqref{eqB}, leads to the inequality for $x\ge x_0$:
$$
\re p \le |p|\le C\re p.
$$

So, if the potential has the estimate \eqref{eqB} and $\kappa>0$, the weights are simplified even more: $|q|^{1/2}$, $1/|q|^{1/2}$.

\section{Auxiliary results}

{\noindent\bf Proof of the Lemma \ref{lmNOpen}.}

Since $0\in\mathscr{N}$, $\mathscr{N}$ is not empty. Let us show that there exists the neighborhood of
zero $\Upsilon=\Upsilon(0)\subset\mathscr{N}$. Without loss of generality, this is sufficient to prove that $\mathscr{N}$ is open.

Let $C>0$ be the constant in \eqref{eqA}.

For $x\ge x_0$ and $|\lambda/q|<1$ we set $p_\lambda=(q-\lambda)^{1/2}=p+O(\lambda/p)$.

From here, all $O$ are considered over a base $\lambda/p\to0$. Since \eqref{eqA} implies $|p|\ge C$ for all $x\ge x_0$, the limit
$\lambda/p\to0$ holds either when $\lambda\to0$ uniformly in $x\ge x_0$, or when $|q|\to+\infty$ as $x\to+\infty$ uniformly in compact
sets in $\lambda$.

We estimate
$$
\Bigl|
\frac{p_\lambda'}{p_\lambda}
\Bigr|=\frac{1}{2}
\Bigl|
\frac{q'}{q-\lambda}
\Bigr|\le
\frac{1}{2}
\Bigl|
\frac{q'}{q}
\Bigl|
(1+O(\lambda/p^2))=
\Bigl|
\frac{p'}{p}
\Bigl|
(1+O(\lambda/p^2)).
$$

It follows from \eqref{eqA} that $|p'/p|=O(p)$ for a.e. $x\ge x_0$, thus taking into account previous estimates:
\begin{equation}
\label{eqplambdeq}
\Bigl|
\frac{p_\lambda'}{p_\lambda}
\Bigr|=
\Bigl|
\frac{p'}{p}
\Bigr|+O(\lambda/p),\quad
p_\lambda=p+O(\lambda/p).
\end{equation}

We take the neighborhood of zero $\Upsilon$ small enough that for all $\lambda\in\Upsilon$ each of $O(\lambda/p)$ in \eqref{eqplambdeq}
does not exceed $C/3$ in absolute value for a.e. $x\ge x_0$.

Then for $\lambda\in\Upsilon$, and a.e. $x\ge x_0$:
\begin{equation}
\label{eqplambdeq2}
\re p_\lambda\ge\frac{C}{2}+\frac{1}{2}\Bigl|
\frac{p_\lambda'}{p_\lambda}
\Bigr|,
\end{equation}
that is, for all $\lambda\in\Upsilon$, the potential $q_\lambda=q-\lambda$ satisfies the condition {\bf A} for $x\ge x_0$
with a common constant $C'=C/2$, hence $\Upsilon\subset\mathscr{N}$.

Now we turn to the second assertion of the Lemma and show the sufficiency.

Let $|q|\to+\infty$ as $x\to+\infty$. We fix an arbitrary $\lambda\in\CC$, and turn to formulas \eqref{eqplambdeq}.

Taking $x_1\ge x_0$ large enough that each of $O(\lambda/p)$
does not exceed $C/3$ in absolute value for a.e. $x\ge x_1$, the required estimate \eqref{eqplambdeq2} is obtained
immediately, which means $\lambda\in\mathcal{N}$.

Now suppose that $\mathcal{N}=\CC$ but $|q|\not\to+\infty$. We take the sequence $x_n\to+\infty$ such that there are finite limits of
$|q(x_n)|$ and $\arg q(x_n)$ as $n\to\infty$. Let $\lambda_0=\lim q(x_n)$. Since $|q(x_n)-\lambda_0|\to0$, then $\re p_{\lambda_0}(x_n)\to0$,
$\lambda_0\not\in\mathcal{N}$. We came to a contradiction.
\qquad$\Box$

\begin{Lemma}
\label{lmst1}
Let $q\in AC_{loc}(\RR_+)$, $|q(x)|>C_1>0$, and $|q(x)'/q^{3/2}(x)|< C_2$ for a.e. $x\ge0$. Then for any solution $y$ of the equation
\begin{equation}
\label{eqmainequat}
-y''+qy=0,\quad x\in\RR_+
\end{equation}
the condition $y\in L_2(\RR_+)$ implies $y'/|q|^{1/2}\in L_2(\RR_+)$.
\end{Lemma}
{\noindent\bf Proof of the Lemma \ref{lmst1}.} Consider $y=y(x)$ --- non-trivial solution to \eqref{eqmainequat}, let
$y\in L_2(\RR_+)$.

We divide both sides of \eqref{eqmainequat} by $|q(x)|$, multiply by $\overline y$, and integrate over $[0,T]$ for arbitrary $T>0$. Taking the first integral by parts,
denoting by $O(1)$ quantities bounded in $T\ge 0$, we obtain the equality:
\begin{equation}
\label{eqspclmst1}
-\frac{\overline y y'}{|q|}\vint{T}{}+\int\limits_0^T y'\Bigl\{
\frac{\overline y'}{|q|}-\frac{|q|'}{|q|^2}\overline y
\Bigr\}\,dx+O(1)=0.
\end{equation}

Let
$$
H(T)=\int\limits_0^T \frac{|y'|^2}{|q|}\,dx.
$$

Taking into account the inequalities: $|q(x)|'\le |q'(x)|$ and $|q'(x)|/|q(x)|^2<C_2/|q(x)|^{1/2}$, we estimate using the Cauchy--Bunyakovsky inequality:
$$
\Bigl|
\int\limits_0^T \frac{|q|'}{|q|^2}\overline y y'\,dx\Bigr|\le C_2\Bigl(\int\limits_0^T|y|^2\,dx
\Bigr)^{1/2}\Bigl(\int\limits_0^T\frac{|y'|^2}{|q|}\,dx
\Bigr)^{1/2}=O(H^{1/2}(T)).
$$

Assume that $H(T)\to+\infty$ for $T\to+\infty$. Then, considering the real parts of both sides of \eqref{eqspclmst1}, we get:
$$
\re\Bigl(\frac{\overline y y'}{|q|}\Bigr)\vint{T}{}=H(T)+O(H^{1/2}(T))+O(1),
$$
whence $\re(\overline y(T) y'(T))\to+\infty$, which is impossible, as $y\in L_2(\RR_+)$, and
$$
\re(\overline y(T) y'(T))=\frac{1}{2}\frac{d}{dx}|y|^2\vint{T}{}.
$$

Thus, $H(T)$ is bounded.\qquad$\Box$
\bigskip

Denote $q_\lambda(x)=q(x)-\lambda$. If $q_\lambda$ satisfies the condition {\bf A} for $x\ge 0$, let $p_\lambda=p_\lambda(x)=\sqrt{q_\lambda(x)}$, choosing the branch so that
for a.e. $x\ge 0$
$$
\re p_\lambda\ge C+\frac{1}{2}\Bigl|\frac{p_\lambda'}{p_\lambda}\Bigr|,
$$
denote
$$
\rho_\lambda=\rho_\lambda(x)=\re p_\lambda(x)-\frac{1}{2}\Bigl|\frac{p_\lambda'(x)}{p_\lambda(x)}\Bigr|.
$$

The following Lemma is directly related to the Weyl--Titchmarsh theory and the construction of the Weyl function $\mu(\lambda)$.

\begin{Lemma}
\label{lm3}
Let for all $\lambda\in\Upsilon$ the function $q_\lambda(x)=q(x)-\lambda$ satisfies the condition {\bf A} for $x\ge x_0$ with a common constant $C>0$,
where $\Upsilon\subset\CC$ is an open and connected set.

Consider the solutions $\varphi(x,\lambda)$, $\psi(x,\lambda)$ to the equation \eqref{eqodnlambda} given by the initial conditions
$$
\varphi(0,\lambda)=0,\quad \varphi'(0,\lambda)=1,\quad \psi(0,\lambda)=1,\quad \psi'(0,\lambda)= 0.
$$

There exists a single-valued function $\mu(\lambda)$ analytic in $\Upsilon$ such that
\begin{itemize}
\item the function $\mu(\lambda)\ne 0$ for all $\lambda\in\Upsilon$;
\item for any $\lambda\in\Upsilon$ the function $w(x,\lambda)=\varphi(x,\lambda)+\mu(\lambda)\psi(x,\lambda)\in L_2(\RR_+)$ as a function of $x$,
and is the only non-trivial solution to the equation
\eqref{eqodnlambda} in $L_2(\RR_+)$ up to a constant factor (not depending on $x$, but possibly depending on $\lambda$);
\item for any $\lambda\in\Upsilon$ the functions $w(x,\lambda)\in L_2(\RR_+,\rho_\lambda)$, $w'(x,\lambda)\in L_2(\RR_+,\rho_\lambda/|q_\lambda|)$ as
the functions of $x$;
\item for $x\to+\infty$ there exists $\lim \re\bigl(\overline{w(x,\lambda)}w'(x,\lambda)/p_\lambda(x)\bigr)=0$.
\end{itemize}
\end{Lemma}
{\noindent\bf Proof of the Lemma \ref{lm3}.} Since the solutions of the Cauchy problem \eqref{eqodnlambda} depend analytically on the parameter
$\lambda$ (see, for example, \cite{Olver}), the solutions
$\varphi(x,\lambda)$, $\psi(x,\lambda)$ --- are continuous in both arguments and entire functions of the parameter $\lambda\in\CC$ for a fixed $x\ge0$.

Without loss of generality, the closure of $\overline{\Upsilon}$ is a compact, and the condition {\bf A} is satisfied for $q_\lambda$ in some
wider open and connected set $O\supset\overline{\Upsilon}$. Otherwise, we will cover
$\Upsilon$ by bounded neighborhoods $\{O_\alpha\}$, each of which satisfies the conditions of the Lemma. Next, we construct own $\mu_\alpha(\lambda)$ for each $O_\alpha$.
At the intersection of such neighborhoods, the functions $\mu_\alpha(\lambda)$ will
coincide in view of the uniqueness of $w(x,\lambda)$ in the class of solutions from $L_2(\RR_+)$ up to a constant factor, as well as the independence of solutions $\varphi$ and
$\psi$.

Taking into account the smoothness of $q\in AC_{loc}(\RR_+)$, we introduce the auxiliary system of equations that is
equivalent to the equation \eqref{eqodnlambda} for $y_1\equiv y$, $y_2\equiv y'/p_\lambda$, $\lambda\in\overline{\Upsilon}$, a.e. $x\ge0$,:
\begin{equation}
\label{eqspecsystem}
\begin{pmatrix}
y_1\\y_2
\end{pmatrix}'=
\begin{pmatrix}
0 & \phantom{p_\lambda'}p_\lambda\\
p_\lambda & -p_\lambda'/p_\lambda
\end{pmatrix}
\begin{pmatrix}
y_1\\y_2
\end{pmatrix}.
\end{equation}
The solutions of the system are understood in the same sense as \cite{NaimarkBk}[Ch.V, \S16, Th. 1].

Multiplying the first equality of the system by $\overline{y}_2$, multiplying the second equality by $\overline{y}_1$,
adding them to each other and integrating over an arbitrary interval $[0,a]$,
we get:
$$
\int\limits_{0}^a
\bigl(y_1'\overline{y}_2+y_2'\overline{y}_1\bigr)\,dx=
\int\limits_{0}^a
\bigl(p_\lambda\,(|y_1|^2+|y_2|^2)-\frac{p_\lambda'}{p_\lambda}\,y_2\overline{y}_1\bigr)\,dx ,
$$
which, after transforming the first integral and separating the real parts, converts to:
\begin{equation}
\label{eqpformqualty}
\re y_1\overline{y}_2\vint{0}{a}=\int\limits_{0}^a
\bigl(\re p_\lambda\,(|y_1|^2+|y_2|^2)-\re(\frac{p_\lambda'}{p_\lambda}\,y_2\overline{y}_1) \bigr)\,dx.
\end{equation}

Applying the estimate
$$
\re(\frac{p_\lambda'}{p_\lambda}\,y_2\overline{y}_1)\le\frac{1}{2}\,\Bigl|\frac{p_\lambda'}{p_\lambda}\Bigr|\,(|y_1|^2+|y_2|^2),
$$
we obtain from \eqref{eqpformqualty}:
\begin{equation}
\label{eqpformineq}
\re y_1\overline{y}_2\vint{0}{a}\ge\int\limits_{0}^a
\rho_\lambda(x)(|y_1|^2+|y_2|^2)\,dx.
\end{equation}

Consider the solutions $U=U(x,\lambda)=(u_1,u_2)^T$ and $V=V(x,\lambda)=(v_1,v_2)^T$ to \eqref{eqspecsystem},
given by the initial conditions:
$$
U(0,\lambda)=\begin{pmatrix}
0\\1
\end{pmatrix},\quad
V(0,\lambda)=\begin{pmatrix}
1\\0
\end{pmatrix},
$$
--- analytic functions in $\lambda\in\overline{\Upsilon}$ for all $x\ge 0$.

Consider the special solution $Y=(y_1,y_2)^T=U+\theta_aV$, where the parameter $\theta_a\in\CC$
we choose from the condition
$\re y_1\overline{y}_2\vint{a}{}=0$, which will be satisfied if one of the two equalities is satisfied:
\begin{equation}
\label{eq1fortheta}
u_2(a)+\theta_av_2(a)=0,
\end{equation}
or
\begin{equation}
\label{eq2fortheta}
\re\frac{u_1(a)+\theta_av_1(a)}{u_2(a)+\theta_av_2(a)}=0.
\end{equation}

Hereinafter we shorten the notation, meaning by $u_j(x)=u_j(x,\lambda)$, $v_j(x)=v_j(x,\lambda)$, $j=1,2$; $\theta_a =\theta_a(\lambda)$.

Since the determinant (aka the Wronskian to \eqref{eqspecsystem})
$$
\begin{vmatrix}
u_1(a) & v_1(a)\\
u_2(a) & v_2(a)
\end{vmatrix}=\frac{p_\lambda(0)}{p_\lambda(a)}\begin{vmatrix}
u_1(0) & v_1(0)\\
u_2(0) & v_2(0)
\end{vmatrix}=-\frac{p_\lambda(0)}{p_\lambda(a)}\ne0,
$$
for fixed $u_1(a)$ and $v_1(a)$ the equality \eqref{eq2fortheta} corresponds to the generalized circle
$\Theta_a=\Theta_a(\lambda)\ni\theta_a(\lambda)$
in the complex $\theta_a$-plane exept for one point $\theta_a$ that is determined by \eqref{eq1fortheta}. Completing $\Theta_a$ with this point,
we simultaneously cover both equalities
\eqref{eq1fortheta} and \eqref{eq2fortheta}.

Let us find out the parameters of the circle $\Theta_a$. Let for a fixed $a>0$
$$
l_a=\frac{u_1(a)+\theta_av_1(a)}{u_2(a)+\theta_av_2(a)},\quad
\theta_a=\frac{u_1(a)-l_au_2(a)}{l_av_2(a)-v_1(a)}.
$$

The point $\theta_a^\infty=\infty$ corresponds to $l_a^\infty=v_1(a)/v_2(a)$. The criterion for the boundedness of the generalized circle $\Theta_a$ is expressed by
the equality $\re v_1\overline{v}_2\vint{a}{}\ne0$. Referring to \eqref{eqpformineq}, taking into account that $V\not\equiv0$ for $x\in[0,a]$,
\begin{equation}
\label{eqrepsi12}
\re v_1\overline{v}_2\vint{a}{}\ge \int\limits_{0}^a
\rho_\lambda(x)(|v_1|^2+|v_2|^2)\,dx\ge C\int\limits_{0}^a
(|v_1|^2+|v_2|^2)\,dx>0,\\
\end{equation}

Therefore, $\Theta_a$ has a finite radius. Since $\re l_a^\infty>0$, the interior of the circle corresponds to
the set of points $\theta_a$ for which the corresponding $l_a$ lie
in the left half-plane. The center of the circle $\theta_a^*$ corresponds to the value
$$
l_a^*=-\overline{v_1(a)}/\overline{v_2(a)},\quad
\theta_a^*=-\frac{u_1(a)\overline{v_2(a)}+\overline{v_1(a)}u_2(a)}{\overline{v_1(a)}v_2(a)+\overline{v_2(a)}v_1(a)},
$$
thanks to \eqref{eqrepsi12}, it's easy to calculate the radius of the circle:
$$
R_a(\lambda)=\Bigl|
-\frac{u_1(a)}{v_1(a)}+\frac{u_1(a)\overline{v_2(a)}+\overline{v_1(a)}u_2(a)}{\overline{v_1 (a)}v_2(a)+\overline{v_2(a)}v_1(a)}
\Bigr|=\frac{1}{2\re v_1\overline{v}_2\vint{a}{}}\le
\frac{C}{\int\limits_{0}^a\rho_\lambda(x)(|v_1|^2+|v_2|^2)\,dx},
$$

Here we note that uniformly in $\lambda\in\overline\Upsilon$ for $a\to+\infty$, each of the integrals
$$
\int\limits_{0}^a\rho_\lambda(x)(|u_1|^2+|u_2|^2)\,dx\to+\infty,\quad\int\limits_{0}^a\rho_ \lambda(x)(|v_1|^2+|v_2|^2)\,dx\to+\infty.
$$

Indeed, for $\lambda\in\overline\Upsilon$, $u_1,v_1\not\in L_2(\RR_+)$, otherwise, in view of Lemma \ref{lmst1}, $u_2=u_1'/p_\lambda\in L_2(\RR_+)$,
(or $v_2=v_1'/p_\lambda\in L_2(\RR_+)$), and then $h_1(x)=\re u_1\overline{u}_2\vint{0}{x}$, (or
$h_2(x)=\re v_1\overline{v}_2\vint{0}{x}$) is integrable over $\RR_+$. If that were the case, then
\eqref{eqpformineq} implies $U\equiv0$ or $V\equiv0$.

Let us demonstrate the uniformity of the limit on the example of the second integral. Under the assumption that the limit is not uniform,
there exist the monotone sequence $a_n\to+\infty$
and the sequence $\lambda_n\to\lambda_0\in\overline\Upsilon$ such that for any $n\in\NN$, $k\ge n$
$$
\int\limits_{0}^{a_n}\rho_{\lambda_k}(x)(|v_1(x,\lambda_k)|^2+|v_2(x,\lambda_k)|^2)\,dx<C.
$$
Taking into account the continuity of $v_j$ in both arguments, we take the limit inside the integral. Then for all $n\in\NN$:
$$
\int\limits_{0}^{a_n}\rho_{\lambda_0}(x)(|v_1(x,\lambda_0)|^2+|v_2(x,\lambda_0)|^2)\,dx\le C,
$$
whence $v_1(x,\lambda_0)\in L_2(\RR_+)$, which, as already noted, is impossible.

Thus, $R_a(\lambda)\to 0$ as $a\to+\infty$ uniformly in $\lambda\in\overline\Upsilon$.

Denote by $K_a=K_a(\lambda)$ the closed disk bounded by $\Theta_a=\partial K_a$. Let us show that for $a_1<a_2$ the inclusion $K_{a_2}\subset K_{a_1}$ takes place.
Considering
$\theta\in\Theta_{a_2}$, we show that for the solution $Y=U+\theta V$ the quantity $\re y_1\overline{y}_2\vint{a_1}{} \le0$.
Using \eqref{eqpformqualty} and the condition $\re y_1\overline{y}_2\vint{a_2}{}=0$, that is equal to $\theta\in\Theta_{a_2}$, we evaluate:
$$
\re y_1\overline{y}_2\vint{a_1}{}=-
\int\limits_{a_1}^{a_2}\bigl(\re p_\lambda\,(|y_1|^2+|y_2|^2)-\re(\frac{p_\lambda'}{p_\lambda }\,y_2\overline{y}_1)\bigr)\,dx\le
-C\int\limits_{a_1}^{a_2}(|y_1|^2+|y_2|^2)\,dx\le0.
$$

Thus, the closed disks $K_a$ are embedded into each other with grows of $a$. Since $R_a\to0$, there us a unique common point $\theta(\lambda)\in\cap K_a(\lambda)$.
We proved, that $\theta_a^0(\lambda)\to\theta(\lambda)$ as
$a\to+\infty$. The quantity $|\theta_a^0(\lambda)-\theta(\lambda)|\le 2R_a(\lambda)\to0$ uniformly in $\lambda\in\overline{\Upsilon}$.
Since each of $\theta_a^0(\lambda)$ is analytic in $\overline{\Upsilon}$, the limit function $\theta(\lambda)$ is also analytic in $\overline{\Upsilon}$.

Let $Y(x,\lambda)=U(x,\lambda)+\theta(\lambda)V(x,\lambda)$. For all $a>0$ the quantity $\re y_1\overline{y}_2\vint{a}{}\le0$. With
\eqref{eqpformineq}, for any $a>0$:
$$
C\int\limits_{0}^{a}(|y_1|^2+|y_2|^2)\,dx\le\int\limits_{0}^{a}\rho_\lambda(x)(| y_1|^2+|y_2|^2)\,dx\le-\re\theta,
$$
taking into account the arbitrariness of $a>0$,
\begin{equation}
\label{eqyinL2}
C\int\limits_{0}^{+\infty}(|y_1|^2+|y_2|^2)\,dx\le\int\limits_{0}^{+\infty}\rho_\lambda( x)(|y_1|^2+|y_2|^2)\,dx\le-\re\theta,
\end{equation}

The constructed solution $Y$ of the system corresponds to the solution $y_1=u_1+\theta v_1$ of the homogeneous equation \eqref{eqodnlambda}.

Since $u_1(x,\lambda)\equiv p_\lambda(0)\varphi(x,\lambda)$ and $v_1(x,\lambda)\equiv \psi(x,\lambda)$, then
$\mu(\lambda)=\theta(\lambda)/p_\lambda(0)$, and $w(x,\lambda)=y_1(x,\lambda)/p_\lambda(0)$.

From \eqref{eqyinL2} we conclude that $w(x,\lambda)\in L_2(\RR_+,\rho_\lambda)$ and $w'(x,\lambda)\in L_2(\RR_+, \rho_\lambda/|q_\lambda|)$ as the functions of $x$.

Lemma \ref{lmst1} implies that the function $h=\re\bigl(\overline{w}w'/p_\lambda\bigr)$ is integrable, and the equality \eqref{eqpformqualty} implies that
there exists $\lim h(x)$ as $x\to+\infty$. Thus $h(x)\to0$.

The fact that $\mu(\lambda)\ne0$ for all $\lambda\in\Upsilon$ follows from the fact that $u_1\not\in L_2(\RR_+)$ for all $\lambda\in\Upsilon$.

The fact that $v_1\not\in L_2(\RR_+)$ implies the uniqueness of $w(x,\lambda)$ in the class of solutions
from $L_2(\RR_+)$ up to a constant factor.\qquad$\Box$

\begin{Lemma}
\label{lmResolvent}
Let $q$ satisfies the condition {\bf A} for $x\ge 0$,

Consider the solutions $\varphi(x)$, $\psi(x)$ to the equation \eqref{eqmainequat} given by the initial conditions
$$
\varphi(0)=0,\quad \varphi'(0)=1,\quad \psi(0)=1,\quad \psi'(0)=0.
$$

Let $w=w(x)=w(x,0)\in L_2(\RR_+)$ be the solution to \eqref{eqmainequat} delivered by Lemma \ref{lm3}.

For $f\in L_2(\RR_+,1/|q|)$ we denote
\begin{gather*}
(R_0f)(x)=w(x)\frac{1}{\mathscr{W}_0}\int\limits_0^x
\varphi(\xi)f(\xi)\,d\xi+
\varphi(x)\frac{1}{\mathscr{W}_0}\int\limits_x^{+\infty}
w(\xi)f(\xi)\,d\xi,\\
(R_1f)(x)=w(x)\frac{1}{\mathscr{W}_1}\int\limits_0^x
\psi(\xi)f(\xi)\,d\xi+
\psi(x)\frac{1}{\mathscr{W}_1}\int\limits_x^{+\infty}
w(\xi)f(\xi)\,d\xi,\\
(R_0'f)(x)=\frac{d}{dx}(R_0f)(x),\quad (R_1'f)(x)=\frac{d}{dx}(R_1f)(x),
\end{gather*}
where $\mathscr{W}_{0,1}$ are the Wronskians of solutions to the homogeneous equation \eqref{eqmainequat}:
$$
\mathscr{W}_0=\begin{vmatrix}
w(\xi) & \varphi(\xi) \\
w'(\xi) & \varphi'(\xi)
\end{vmatrix}\ne0,\quad
\mathscr{W}_1=\begin{vmatrix}
w(\xi) & \psi(\xi) \\
w'(\xi) & \psi'(\xi)
\end{vmatrix}\ne0.
$$

Then $R_{0,1}$ are bounded operators from $L_2(\RR_+,1/|q|)$ to $L_2(\RR_+,\rho)$,
$R_{0,1}'$ --- bounded operators from $L_2(\RR_+,1/|q|)$ to $L_2(\RR_+,\rho/|q|)$.

Operators $L_{U_0}$ with Dirichlet boundary conditions $U_0(y)=y(0)$ and $L_{U_1}$ with Neumann boundary conditions $U_1(y)=y'(0)$ are bounded invertible in $L_2 (\RR_+)$,
$L_{U_{0,1}}^{-1}=R_{0,1}\vint{L_2(\RR_+)}{}$.
\end{Lemma}
{\noindent\bf Proof of the Lemma \ref{lmResolvent}.} It is shown in the proof of the Lemma \ref{lm3}, that
if $q$ satisfies the condition {\bf A} for $x\ge 0$, then $\varphi,\psi\not\in L_2(\RR_+)$.

Thus, each of the pairs of solutions $\{w,\varphi\}$, $\{w,\psi\}$ is linearly independent, and the corresponding Wronskians are $\mathscr{W}_{0,1}\ne0$.

Further proof will be carried out only for the operator $R_1$ due to the verbatim repetition in the case of $R_0$. We will omit the index 1 in
the definition of operators to simplify the notation.

Let $f\in L_2(\RR_+,1/|q|)$,  consider the inhomogeneous equation
\begin{equation}
\label{eqmainequatnood}
-y''+qy=f.
\end{equation}

It is directly verified that $y,y'\in AC_{loc}(\RR_+)$ for $y=Rf$, and $y$ is the solution to \eqref{eqmainequatnood}.

Like in the proof of the Lemma \ref{lm3}, we associate with \eqref{eqmainequatnood} the equivalent auxiliary inhomogeneous system of the first-order equations
for $y_1\equiv y$, $y_2\equiv y'/p$:
$$
\begin{pmatrix}
y_1\\y_2
\end{pmatrix}'-
\begin{pmatrix}
0 & \phantom{p'}p\\
p & -p'/p
\end{pmatrix}
\begin{pmatrix}
y_1\\y_2
\end{pmatrix}+
\begin{pmatrix}
0\\g
\end{pmatrix}=0,
$$
where $g=f/p\in L_2(\RR_+)$.

Let us transform the inhomogeneous system
similar to the transformation \eqref{eqpformqualty} of the homogeneous system:
\begin{equation}
\label{eqpformineq2}
\re y_1\overline{y}_2\vint{0}{a}=\int\limits_0^a
\bigl(\re p\,(|y_1|^2+|y_2|^2)-\re(\frac{p'}{p}\,y_2\overline{y}_1)\bigr)\,dx
-\int\limits_0^a\re\bigl(g\overline{y}_1\bigr)\,dx
\end{equation}

By $f^N$ we denote the cut: $f^N(x)=f(x)$ for $x\in[0,N]$, $f^N(x)=0$ for $x>N$. Denote $g^N=f^N/p$, $y^N=Rf^N$, $y_1^N\equiv y^N$,
$y_2^N\equiv(y^N)'/p$.

Since for $x>N$:
$$
y^N=w(x)A^N,\quad A^N=\frac{1}{\mathscr{W}}\int\limits_0^N
\psi(\xi)f^N(\xi)\,d\xi,
$$
then $y^N\in L_2(\RR_+,\rho)$, $\re y^N_1\overline{y}^N_2\vint{a}{}\to0$ as $a\to+\infty$. The latter follows from Lemma \ref{lm3}.

Applying \eqref{eqpformineq2} to $y_{1,2}^N$, taking into account that $\re y^N_1\overline{y}^N_2\vint{0}{}=0$, we conclude
that all integrals converge on $+\infty$, and
$$
\int\limits_0^{+\infty}
\bigl(\re p\,(|y_1^N|^2+|y_2^N|^2)-\re(\frac{p'}{p}\,y_2^N\overline{y}_1^ N)\bigr)\,dx\le \|g^N\|\cdot\|y^N\|,
$$
where $\|\cdot\|$ means the norm in $L_2(\RR_+)$. Just like we came from \eqref{eqpformqualty} to \eqref{eqpformineq},
$$
\|y_1^N\|^2_\rho+\|y_2^N\|^2_\rho=\int\limits_0^{+\infty}\rho(x)(|y_1^N|^2+|y_2^ N|^2)\le
\int\limits_0^{+\infty}
\bigl(\re p\,(|y_1^N|^2+|y_2^N|^2)-\re(\frac{p'}{p}\,y_2^N\overline{y}_1^ N)\bigr)\,dx,
$$
where $\|\cdot\|_\rho$ means the norm in the weighted $L_2(\RR_+,\rho)$.

As a result, taking into account \eqref{eqA}, we get:
$$
C\|y^N\|^2\le\|y_1^N\|^2_\rho+\|y_2^N\|^2_\rho\le \|g^N\|\cdot\|y^N\|,
$$
where $C$ does not depend on either $g$ or $N$.

For $M,N>0$, considering the function $g^N-g^M$ instead of $g^N$, one can similarly obtain the estimates:
$$
C\|y^N-y^M\|^2\le\|y_1^N-y_1^M\|_\rho^2+\|y_2^N-y_2^M\|_\rho^2\le\|g^N-g^M\|\cdot\|y^N-y^M\|,
$$
proving that $y_{1,2}^N$ is fundamental both in $L_2(\RR_+)$ and $L_2(\RR_+,\rho)$ as $N\to+\infty$. Since
$y_{1,2}^N\Rightarrow y_{1,2}$ as $N\to+\infty$
uniformly in any segment, so
$y_{1,2}^N\to y_{1,2}$ in the sense of each of the norms $L_2(\RR_+)$ and $L_2(\RR_+,\rho)$, and
\begin{equation}
\label{eqlm04spr}
C\|y\|\le\|g\|,\quad \|y_1\|^2_\rho+\|y_2\|^2_\rho\le \|g\|\cdot\|y\|,\quad C^{1/2}\|y_{1,2}\|_\rho\le\|g\|=\|f\|_{1/|q|},
\end{equation}
which proves that the operators $R$ and $R'$ are bounded in the corresponding spaces.

Let us prove the injectivity of $L_U$. If $L_Uw_1=L_Uw_2=f$ then $w_1-w_2=C\psi$, since $\psi\not\in L_2(\RR_+)$ we conclude that $C=0$.

Since $R\vint{L_2(\RR_+)}{}$ is defined in the entire $L_2(\RR_+)$, then
$L_U^{-1}=R\vint{L_2(\RR_+)}{}$.\qquad$\Box$

\section{Proofs of the Theorems}

{\noindent\bf Proof of the Theorem \ref{thWeil}.} Let's split the proof into several steps.

\bigskip
{\bf\noindent Step 1}. Construction of the function $\eta(x,\lambda)$.

It follows from Lemma \ref{lmNOpen} that for each $\lambda_0\in\Omega$ there exist the bounded neighborhood $\Upsilon(\lambda_0)$, $a\ge0$
such that $q_\lambda(x)=q(x)-\lambda$ satisfies the condition {\bf A} for $x\ge a$,
with a common constant $C_1$ for all $\lambda\in\Upsilon(\lambda_0)$.

Based on Lemma \ref{lm3}, we construct the solution $w_a(x,\lambda)=\varphi_a(x,\lambda)+\mu(\lambda)\psi_a(x,\lambda)$,
delivered by this Lemma, $w_a(x,\lambda)\in L_2(a,+\infty)$ for all $\lambda\in\Upsilon(\lambda_0)$ as a function of $x$.

Here
$\varphi_a$, $\psi_a$ --- the solutions to the Cauchy problem given by the conditions
$$
\varphi_a(a,\lambda)=0,\quad \varphi_a'(a,\lambda)=1,\quad \psi_a(a,\lambda)=1,\quad \psi_a'(a,\lambda)= 0,
$$
--- entire functions of the argument $\lambda$ for each $x\in\RR_+$; the function $\mu(\lambda)$ is analytic in $\Upsilon(\lambda_0)$.

The sets of points $\lambda\in\Upsilon(\lambda_0)$ for which $w_a(0,\lambda)=0$ or $w_a'(0,\lambda)=0$ are discrete, being zeros of
analytical functions.

Outside these points, consider two functions:
\begin{equation*}
\begin{split}
\eta_0(x,\lambda)&=\varphi(x,\lambda)+\nu_0(\lambda)\psi(x,\lambda),\quad \nu_0(\lambda)=
\frac{w_a(0,\lambda)}{w_a'(0,\lambda)}\\
\eta_1(x,\lambda)&=\nu_1(\lambda)\varphi(x,\lambda)+\psi(x,\lambda),\quad \nu_1(\lambda)=
\frac{w_a'(0,\lambda)}{w_a(0,\lambda)}.
\end{split}
\end{equation*}

By construction, $\eta_{0,1}\in L_2(\RR_+)$ as functions of $x$. The functions $\nu_{0,1}$ are meromorphic in
$\Upsilon(\lambda_0)$ and are uniquely determined due to the uniqueness of the solution in $L_2$ (Lemma \ref{lm3}).

For the same reason, the functions $\nu_{0,1}$ can be uniquely extended as meromorphic to $\Omega$, where $\nu_0\nu_1\equiv 1$.

By $n_0(\lambda)$
we denote the canonical product constructed by the zeros of $\nu_0$ in $\Omega$ taking into account their multiplicity. Let
$\eta(x,\lambda)=n_0(\lambda)\eta_1(x,\lambda)$.

By construction, $\eta(x,\lambda)$ is the non-trivial solution to \eqref{eqodnlambda} for all $\lambda\in\Omega$; $\eta$, $\eta'$ are continuous in both arguments, analytic
in $\Omega$ as functions of $\lambda$ for each $x\in\RR_+$. For each $\lambda\in\Omega$, the function $\eta\in L_2(\RR_+)$ as a function of $x$.
Uniqueness in $L_2(\RR_+)$ follows from Lemma \ref{lm3}.

\bigskip
{\bf\noindent Step 2}. Construction of the resolvent.

By construction, $\mathscr{W}(\lambda)$ is analytic in $\Omega$, its zeros $\lambda_n\in\Omega$ determine the eigenfunctions $\eta(x,\lambda_n)$ corresponding to
eigenvalues $\{\lambda_n\}$.

We fix $\lambda\not\in\{\lambda_n\}$.

It is easy to check that $R_\lambda f$ is the solution to \eqref{eqmainequatnood}; for any
$f\in L_2(\RR_+)$, $R_\lambda f\in D_U$. Let us show that the operator $R_\lambda$ is bounded.

Since $\lambda\in\mathcal{N}$ then for some $a\ge0$ the function $q_\lambda(x)$ satisfies the condition {\bf A} for $x\ge a$.
According to Lemma \ref{lmResolvent} the following operators are bounded in $L_2(a,+\infty)$:
\begin{gather*}
(R_0f)(x)=\eta(x,\lambda)\frac{1}{\mathscr{W}_0(\lambda)}\int\limits_a^x
\varphi_a(\xi,\lambda)f(\xi)\,d\xi+
\varphi_a(x,\lambda)\frac{1}{\mathscr{W}_0(\lambda)}\int\limits_x^{+\infty}
\eta(\xi,\lambda)f(\xi)\,d\xi,\\
(R_1f)(x)=\eta(x,\lambda)\frac{1}{\mathscr{W}_1(\lambda)}\int\limits_a^x
\psi_a(\xi,\lambda)f(\xi)\,d\xi+
\psi_a(x,\lambda)\frac{1}{\mathscr{W}_1(\lambda)}\int\limits_x^{+\infty}
\eta(\xi,\lambda)f(\xi)\,d\xi,
\end{gather*}
where
$$
\mathscr{W}_0(\lambda)=\begin{vmatrix}
\eta(\xi,\lambda) & \varphi_a(\xi,\lambda) \\
\eta'(\xi,\lambda) & \varphi_a'(\xi,\lambda)
\end{vmatrix}\ne0,\quad
\mathscr{W}_1(\lambda)=\begin{vmatrix}
\eta(\xi,\lambda) & \psi_a(\xi,\lambda) \\
\eta'(\xi,\lambda) & \psi_a'(\xi,\lambda)
\end{vmatrix}\ne0.
$$

Using $R_0$ as an example, we will show that the same formulas define bounded operators in $L_2(\RR_+)$. Let $f\in L_2(\RR_+)$, represent
$f=f_0+f_1$, where $f_0=0$ for $x>a$ and $f_1=0$ for $x\in[0,a]$. Hereinafter, $\|\cdot\|$ will denote the norm in $L_2(\RR_+)$.
\begin{gather*}
\|R_0f_0\|^2=\|R_0f_0\|_{L_2[0,a]}^2\le C\|f_0\|_{L_2[0,a]}^2\le C\|f_0\|^2,\\
\|R_0f_1\|^2\le\Bigl\|\varphi_a(x,\lambda)\frac{1}{\mathscr{W}_0(\lambda)}\int\limits_a^{+\infty}
\eta(\xi,\lambda)f_1(\xi)\,d\xi\Bigr\|_{L_2[0,a]}^2 +\|R_0f_1\|_{L_2[a,+\infty) }^2\le C\|f_1\|^2,
\end{gather*}
the boundedness of the operators in $L_2[0,a]$ is beyond doubt in view of the continuity of the solutions to \eqref{eqodnlambda}
and, consequently, the boundedness on the segment $x\in[0,a]$.

As a result, for all $f\in L_2(\RR_+)$
$$
\|R_0f\|^2=\|R_0f_0\|^2+\|R_0f_1\|^2\le C\|f\|^2,
$$
which proves that $R_0$ is bounded. Similar reasoning proves the boundedness of $R_1$.

For the same reasons, the operators $\tilde R_0$ and $\tilde R_1$ are also bounded:
\begin{gather*}
(\tilde R_0f)(x)=\eta(x,\lambda)\frac{1}{\mathscr{W}_0(\lambda)}\int\limits_0^x
\varphi_a(\xi,\lambda)f(\xi)\,d\xi+
\varphi_a(x,\lambda)\frac{1}{\mathscr{W}_0(\lambda)}\int\limits_x^{+\infty}
\eta(\xi,\lambda)f(\xi)\,d\xi,\\
(\tilde R_1f)(x)=\eta(x,\lambda)\frac{1}{\mathscr{W}_1(\lambda)}\int\limits_0^x
\psi_a(\xi,\lambda)f(\xi)\,d\xi+
\psi_a(x,\lambda)\frac{1}{\mathscr{W}_1(\lambda)}\int\limits_x^{+\infty}
\eta(\xi,\lambda)f(\xi)\,d\xi,
\end{gather*}
and, as a consequence, $R_\lambda$ is bounded, which is their linear combination for $\lambda\not\in\{\lambda_n\}$.

\bigskip
{\bf\noindent Step 3}. Let us show that $\lambda_n$ are the eigenvalues of finite algebraic multiplicity.

Let $\lambda_0\in\{\lambda_n\}$. Since the space of $L_2(\RR_+)$--solutions to the homogeneous equation is one-dimensional, any operator
$L_V$ with boundary conditions defined by some form $V(y)=\beta_0 y(0)+\beta_1 y'(0)$ with linearly independent coefficients $(\beta_0,\beta_1)$
with respect to $(\alpha_0,\alpha_1)$ will not have an eigenvalue at the point $\lambda_0$ and even entirely in some neighborhood of $\Upsilon(\lambda_0)$.

We take two such forms $V_1$ and $V_2$ which are pairwise linearly independent with the form $U$ and between each other.

To simplify the reasoning, we assume that these are $V_1=\alpha_0 y(0)$ and $V_2=\alpha_1 y'(0)$, which occurs when $\alpha_0\alpha_1\ne 0$.
The resolvents of the corresponding operators $R_{V_1,\lambda}$, $R_{V_2,\lambda}$ are analytic in $\Upsilon(\lambda_0)$ operator--functions:
\begin{gather*}
(R_{V_1,\lambda}f)(x)=\eta(x,\lambda)\frac{1}{\mathscr{W}_0(\lambda)}\int\limits_0^x
\alpha_0\varphi(\xi,\lambda)f(\xi)\,d\xi+
\alpha_0\varphi(x,\lambda)\frac{1}{\mathscr{W}_0(\lambda)}\int\limits_x^{+\infty}
\eta(\xi,\lambda)f(\xi)\,d\xi,\\
(R_{V_2,\lambda}f)(x)=\eta(x,\lambda)\frac{1}{\mathscr{W}_1(\lambda)}\int\limits_0^x
\alpha_1\psi(\xi,\lambda)f(\xi)\,d\xi+
\alpha_1\psi(x,\lambda)\frac{1}{\mathscr{W}_1(\lambda)}\int\limits_x^{+\infty}
\eta(\xi,\lambda)f(\xi)\,d\xi,
\end{gather*}
where
$$
\mathscr{W}_0(\lambda)=\begin{vmatrix}
\eta(\xi,\lambda) & \alpha_0\varphi(\xi,\lambda) \\
\eta'(\xi,\lambda) & \alpha_0\varphi'(\xi,\lambda)
\end{vmatrix}\ne0,\quad
\mathscr{W}_1(\lambda)=\begin{vmatrix}
\eta(\xi,\lambda) & \alpha_1\psi(\xi,\lambda) \\
\eta'(\xi,\lambda) & \alpha_1\psi'(\xi,\lambda)
\end{vmatrix}\ne0,
$$
thus
$$
\mathscr{W}(\lambda) R_\lambda =\mathscr{W}_0(\lambda)R_{V_1,\lambda}-\mathscr{W}_1(\lambda)R_{V_2,\lambda}
$$
--- is also analytic in $\lambda\in\Upsilon(\lambda_0)$, hence the poles of $R_\lambda$ coincide with the zeros of $\mathscr{W}(\lambda)$, they have finite multiplicity
that coincides with the dimension of root subspaces \cite[Theor 10.8]{Taylor}.

If $\alpha_0\alpha_1=0$, then if $\alpha_1=0$ we take $V_1=\alpha_0 y(0)+y'(0)$ and $V_2=y'(0)$,
and then $\mathscr{W}(\lambda) R_\lambda =\mathscr{W}_0(\lambda)R_{V_1,\lambda}+\mathscr{W}_1(\lambda)R_{V_2,\lambda }$.
If $\alpha_0=0$ --- we use the forms $V_1=y(0)+\alpha_1 y'(0)$ and $V_2=y(0)$,
then $\mathscr{W}(\lambda) R_\lambda =\mathscr{W}_0(\lambda)R_{V_2,\lambda}-\mathscr{W}_1(\lambda)R_{V_1,\lambda}$.

\bigskip
{\bf\noindent Step 4}. Let us show that if for some $\lambda\in\CC$ and $a\ge0$ the function
$q_\lambda(x)$ satisfies the condition {\bf A} for $x\ge a$, then the definite case is realized.

Since the set of such $\lambda$ is open, without loss of generality
$\lambda$ is not the eigenvalue of $L_U$, i.e. $\eta(x,\lambda)$ does not satisfy the boundary conditions $U$ for $x=0$.

It's sufficient to show that $\overline{L_U^0-\lambda I}=L_U-\lambda I$. Further, without loss of generality, we assume that $\lambda=0$.

The image of the operator $L_U$ is the entire $L_2(\RR_+)$, it's sufficient to show
that the image of $L_U^0$ is dense in $L_2(\RR_+)$. Assuming the contrary, we find $g\in L_2(\RR_+)$ such that for all $y\in D_U^0$
the inner product $(L_U^0y,g)=0$.

Consider the operator $L^0$ --- the restriction of $L_U^0$ to the domain $D^0$ consisting of
finite functions $y\in D_U^0$ for which $y(0)=y'(0)=0$.

In particular, for all $y\in D^0$ the inner product $(L^0y,g)=0$. That means $g$ lies in the domain of the adjoint operator
$(L^0)^*=\overline L$, which, according to \cite{Naimark},
is defined by the expression $\overline l(y)=-y''+\overline{q(x)}y$, where $\overline{q(x)}$ is the complex conjugate potential.
Consequently, $\overline g\in D$ and $l(\overline g)=0$, hence $\overline g=C\eta$ in view of the uniqueness of the solution $\eta\in L_2(\RR_+)$ .

Next, we use the finiteness of the functions $y\in D_U^0$ and the condition $l(\overline g)=0$. Integrating by parts, we obtain for all $y\in D_U^0$:
$$
(L_U^0y,g)=\int\limits_0^\infty (-y''+q(x)y)\overline{g(x)}\,dx =y'(0)\overline g(0) -y(0)\overline g'(0),
$$
hence $\overline g=C\eta$ satisfies the boundary conditions $U$ at zero, hence $C=0$, and $g=0$, since by assumption 0 is not the eigenvalue.
\qquad$\Box$

\bigskip
{\noindent\bf Proof of the Theorem \ref{th02}.} Since the condition $\rho\to+\infty$ implies $|q|\to+\infty$, then according to
Lemma \ref{lmNOpen}, $\mathcal{N}=\CC$.
Without loss of generality, $0\not\in\sigma(L_U)$.

It follows from Theorem \ref{thWeil} that the spectrum of $L_U$ is discrete. This also applies to the operators $L_{x_0,0}$ and $L_{x_0,1}$ in $L_2(x_0,+\infty)$,
given by Dirichlet and Neumann boundary conditions at the point $x_0$.

Consider the operators $R_0$, $R_1$ in $L_2(x_0,+\infty)$ inverse to $L_{x_0,0}$ and $L_{x_0,1}$ respectively.
As noted in the proof of Theorem \ref{thWeil} (Step 2), they extend to bounded operators in $L_2(\RR_+)$.

Taking into account the explicit form of the inverse operator $R=L_U^{-1}$, which follows from Theorem \ref{thWeil}, as well as the explicit forms of $R_0$ and $R_1$, it is easy
to show that for $f\in L_2(\RR_+)$
$$
Rf=c_1 \eta\int\limits_0^{x_0}\chi f\,d\xi+c_2R_0f+c_3R_1f,
$$
where $c_j$ are constants independent of $f$.

The first term in this formula is a one-dimensional operator, hence it is a compact operator. To prove the compactness of $R$
it's sufficient to show the compactness of each of $R_0$, $R_1$ as operators in $L_2(\RR_+)$.

Let's do this using $R_0$ as an example. Consider the bounded set of functions $\mathcal{F}\subset L_2(\RR_+)$, $\|f\|<C$ for all $f\in\mathcal{F}$.
It's sufficient to prove that the set $\mathcal{Y}=\set{y=R_0f}{f\in\mathcal{F}}$ is precompact.

For each $\varepsilon>0$ we take $T\ge x_0$ so that $\rho(x)>1/\varepsilon^2$ for all $x\ge T$. For $f\in\mathcal{F}$ let $y=R_0f$,
then in view of Lemma \ref{lmResolvent}
$$
\int\limits_T^{+\infty}\rho(x)|y(x)|^2\,dx\le \int\limits_{x_0}^{+\infty}\rho(x)|y(x )|^2\,dx\le C_1\int\limits_{x_0}^{+\infty}|f(x)|^2\,dx\le C_1 C,
$$
whence $\|y\|^2_{L_2[T,+\infty)}\le C_1C\varepsilon^2$.

We represent
\begin{equation}
\begin{split}
(R_0f)(x)=&\eta(x)\frac{1}{\mathscr{W}_0}\int\limits_{x_0}^x
\varphi_{x_0}(\xi)f(\xi)\,d\xi-\varphi_{x_0}(x)\frac{1}{\mathscr{W}_0}\int\limits_{x_0}^x\eta(\xi)f(\xi)\,d\xi+\\
+&\varphi_{x_0}(x)\frac{1}{\mathscr{W}_0}\int\limits_0^{+\infty}
\eta(\xi)f(\xi)\,d\xi
=(Kf)(x)+\varphi_{x_0}(x)\frac{1}{\mathscr{W}_0}\int\limits_0^{+\infty}
\eta(\xi)f(\xi)\,d\xi,
\end{split}
\label{eqr0compact}
\end{equation}
where $\eta,\varphi_{x_0}$ are solutions to \eqref{eqmainequat}, $\eta\in L_2(\RR_+)$, $\varphi_{x_0}(x_0)=0 $, $\varphi'_{x_0}(x_0)=1$,
$\mathscr{W}_0=\eta(x_0)$ --- their Wronskian.

The last term in \eqref{eqr0compact} can be considered as a one-dimensional operator from $L_2(\RR_+)$ to $L_2[0,T]$,
and the first --- as a compact operator, acting in the same spaces, which is obvious,
due to the continuity of $\eta$ and $\varphi_{x_0}$ on $[0,T]$. Therefore, the restriction to $[0,T]$ of functions from $\mathcal{Y}$ is precompact in $L_2[0,T]$.
Denote this set of functions as $\mathcal{Y}\vint{[0,T]}{}$.

Let
the restrictions of $\{y_j\vint{[0,T]}{}\}_{j=1}^n$ represent the
finite $\varepsilon\sqrt{C_1C}$-net in $\mathcal{Y}\vint{[0,T]}{}$ in terms of
$L_2[0,T]$ norm. Consider the functions $\{y_j\}$ in $L_2(\RR_+)$.

For each $y\in\mathcal{Y}$ we find the element $y_s\in\{y_j\}$ such that $\|y-y_s\|_{L_2[0,T]}<\varepsilon\sqrt{C_1C }$, then
$$
\|y-y_s\|^2_{L_2(\RR_+)}<C_1C\varepsilon^2+\|y-y_s\|^2_{L_2[T,+\infty)}\le5C_1C\varepsilon^2,
$$
thus $\{y_j\}_{j=1}^n$ form a finite $\varepsilon\sqrt{5C_1C}$-net in $\mathcal{Y}$ in terms of $L_2(\RR_+)$ norm, i.e. $\mathcal{Y}$ is precompact.\qquad$\Box$

\bigskip
{\noindent\bf Proof of the Theorem \ref{th03}.} Without loss of generality, $q(x)\in\Pi_\kappa$ for $x\ge0$. Otherwise
we adjust the potential by adding the constant $C>0$.

We can also find $x_1\ge x_0$ to make \eqref{eqth03kappa} valid for
adjusted potential $q(x)+C$ for $x\ge x_1$, possibly with the new constant $\delta'\in(\delta,1)$. Indeed, for $x\ge x_0$:
$$
\Bigl|\frac{q'}{(q+C)^{3/2}}\Bigr|=\Bigl|\frac{q'}{q^{3/2}}\Bigr|\Bigl| \frac{q}{q+C}\Bigr|^{3/2},
$$
and the last factor tends to 1 as $x\to+\infty$.

So, without loss of generality let $q(x)\in\Pi_\kappa$ for all $x\ge0$, and \eqref{eqth03kappa} be satisfied for all $x\ge x_0$.

Consider $\lambda<0$, denote $q_\lambda=q-\lambda$, $p_\lambda=\sqrt{q_\lambda}$, choosing the branch of the root to satisfy $\re p_\lambda>C>0$ for all $x\ge0$.

Since $\kappa>0$, the following inequalities hold for all $x\ge0$:
$$
|q_\lambda|\ge|q|\tan\kappa,\quad |q_\lambda|\ge|\lambda|\sin\kappa,\quad\re p_\lambda\ge|p_\lambda|\sin\frac{\kappa}{2},
$$
thus
$$
\re p_\lambda\ge|q|^{1/2}\tan^{1/2}\kappa\,\sin\frac{\kappa}{2},
\quad\re p_\lambda\ge|\lambda|^{1/2}\sin^{1/2}\kappa\,\sin\frac{\kappa}{2}.
$$

From these inequalities we conclude that $\re p_\lambda\to+\infty$ either as $x\to+\infty$ uniformly in $\lambda<0$, or
as $\lambda\to-\infty$ uniformly in $x\ge0$.

Let's estimate
$$
\Bigl|\frac{p'_\lambda}{p_\lambda}\Bigr|=\frac{1}{2}\Bigl|\frac{q'_\lambda}{q_\lambda}\Bigr|=
\frac{1}{2}\Bigl|\frac{q'}{q}\Bigr|\Bigl|\frac{q}{q_\lambda}\Bigr|\le\frac{1}{2}\Bigl|\frac{q'}{q}\Bigr|\frac{1}{\tan\kappa}.
$$
This estimation does not depend on $\lambda$, and since $q\in C^1(\RR_+)$, the quantity on the right hand side is bounded at $x\in[0,x_0]$.

Using the above inequalities,
$$
\rho_\lambda=\re p_\lambda-\frac{1}{2}\Bigl|\frac{p'_\lambda}{p_\lambda}\Bigr|=\re p_\lambda\Bigl(
1-\frac{1}{\re p_\lambda}\frac{1}{2}\Bigl|\frac{p'_\lambda}{p_\lambda}\Bigr|
\Bigr).
$$

We take $\lambda_0<0$ so large in absolute value that for all $\lambda\le\lambda_0$ and $x\in[0,x_0]$ the value $\rho_\lambda>C_1\re p_\lambda\ge C|\lambda|^{1/2}$ .

Further for $x\ge x_0$:
\begin{gather}
\notag
\rho_\lambda\ge\re p_\lambda\Bigl(1-\frac{1}{2}\Bigl|\frac{p'_\lambda}{p_\lambda}\Bigr|\frac{1}{|q|^{1/2}\tan^{1/2}\kappa\,\sin(\kappa/2)}\Bigr)\ge\\
\ge\re p_\lambda\Bigl(1-\Bigl|\frac{q'}{q^{3/2}}\Bigr|\frac{1}{4\tan^{3/2}\kappa\,\sin(\kappa/2)}\Bigr)\ge\re p_\lambda(1-\delta),
\label{eqprvth3evdelt}
\end{gather}
whence it follows that $\rho_\lambda\to+\infty$ as $x\to+\infty$, and we can apply Theorem \ref{th03}.

The choice of $\lambda_0<0$ and \eqref{eqprvth3evdelt} imply that for $\lambda\le\lambda_0$ the potential
$q_\lambda(x)$ satisfies the condition {\bf A} for $x\ge 0$, moreover, $\rho_\lambda>C|\lambda|^{1/2}$.

For $\lambda\in\CC$ we denote by $R_{0,\lambda}$ the resolvent of the operator $L_{U_0}$ with Dirichlet boundary conditions and consider its value on the function $f$,
that appears in the statement of the Theorem.
Let $y=R_{0,\lambda}f$. Since the integral \eqref{eqth03int} is equal to zero, the explicit form of the formula for $y$ is simplified:
$$
y(x,\lambda)=\eta(x,\lambda)\frac{1}{\mathscr{W}_0}\int\limits_0^x
\varphi(\xi,\lambda)f(\xi)\,d\xi-
\varphi(x,\lambda)\frac{1}{\mathscr{W}_0}\int\limits_0^x
\eta(\xi,\lambda)f(\xi)\,d\xi,
$$
where $\varphi$ is the solution to \eqref{eqmainequat}, $\varphi(0,\lambda)=0$, $\varphi'(0,\lambda)=1$ and
$\mathscr{W}_0=\eta(0,\lambda)$ is the Wronskian of $\eta$ and $\varphi$.

Obviously, $y(x,\lambda)$ is the solution to the Cauchy problem for $-y''+(q-\lambda)y=f$ with zero initial conditions both
for the function and for derivative.
In view of the uniqueness of such solution, it coincides with
$$
y(x,\lambda)=\psi(x,\lambda)\int\limits_0^x
\varphi(\xi,\lambda)f(\xi)\,d\xi-
\varphi(x,\lambda)\int\limits_0^x
\psi(\xi,\lambda)f(\xi)\,d\xi,
$$
where $\psi$ is the solution to \eqref{eqmainequat}, $\psi(0,\lambda)=1$, $\psi'(0,\lambda)=0$.

For a fixed $x\ge0$, each of the functions $\varphi$, $\psi$ is an entire function of the argument $\lambda$ and has the growth order at most 1/2 \cite{NaimarkBk}.
Moreover, this estimate of the growth order is uniform in segments in $x$.
The function $y$ has the same property, as well as the function
$$
\upsilon(x,\lambda)=\int\limits_0^x y(\xi,\lambda)\,d\xi.
$$

Let us turn to the second formula of \eqref{eqlm04spr}, which we write for $\lambda<-\lambda_0$ as follows:
$$
\|y\|_{\rho_\lambda}^2\le\|f\|_{1/|q_\lambda|}\|y\|.
$$
Applying the estimates $\rho_\lambda>C_1|\lambda|^{1/2}$, $|q_\lambda|\ge C_2|\lambda|$, valid for $x\ge0$, we obtain
$$
C_1|\lambda|^{1/2}\|y\|^2\le \frac{\|f\|}{C_2^{1/2}|\lambda|^{1/2}}\| y\|,
$$
whence $\|y\|=O(1/\lambda)$ as $\RR\ni\lambda\to-\infty$.

Since
$$
|\upsilon(x,\lambda)|\le\sqrt{x}\,\|y\|,
$$
the function $\upsilon(x,\lambda)=O(1/\lambda)$ for fixed $x$ as $\RR\ni\lambda\to-\infty$, at the same time, this entire function has the growth order
not higher than 1/2,
from the Phragmen--Lindelöf theorem we conclude, that $\upsilon(x,\lambda)\equiv 0$, hence $y(\xi,\lambda)\equiv0$, and hence also $f\equiv 0$.\qquad $\Box$

\bigskip
{\noindent\bf Proof of the Theorem \ref{th04}.} Follows directly from the proof of the Lemma \ref{lmResolvent}.\qquad$\Box$
\bigskip

This research was supported by RSF grant No 20-11-20261.

\end{document}